\documentclass[a4paper,12pt]{amsart}
\usepackage{amssymb,amscd,amsmath,a4wide,graphicx,stmaryrd,fullpage,setspace,microtype,xr,ulem,textcomp,enumitem}
\usepackage[pagebackref=true]{hyperref}
\usepackage[all]{xy}
\usepackage[T1]{fontenc}
\usepackage{lmodern}

\title{Six model categories for directed homotopy}

\author[P. Gaucher]{Philippe Gaucher}

\address{Universit\'e de Paris, CNRS, IRIF, F-75006, Paris, France}

\urladdr{http://www.irif.fr/{\~{}}gaucher} 



\swapnumbers

\newcommand{\C}{\mathcal{C}}

\newcommand{\K}{\mathcal{K}}

\newcommand{\W}{\mathcal{W}}
\newcommand{\F}{\mathcal{F}}

\newcommand{\p}{\times}
\renewcommand{\vec}{\overrightarrow}
\renewcommand{\P}{\mathbb{P}}

\newtheorem*{thmN}{Theorem}

\newtheorem{thm}{Theorem}[section]
\newtheorem{prop}[thm]{Proposition}

\newtheorem{cor}[thm]{Corollary}

\newtheorem{defn}[thm]{Definition}
\newtheorem{nota}[thm]{Notation}
\newtheorem{rem}[thm]{Remark}
\newcommand{\bd}{\begin{defn}}
	\newcommand{\ed}{\end{defn}}
\newcommand{\bp}{\begin{prop}}
	\newcommand{\ep}{\end{prop}}
\newcommand{\bth}{\begin{thm}}
	\renewcommand{\eth}{\end{thm}}
\newcommand{\bpf}{\begin{proof}}
	\newcommand{\epf}{\end{proof}}
\newcommand{\bc}{\begin{cor}}
	\newcommand{\ec}{\end{cor}}

\newcommand{\fL}[1]{\ar@{->}[ll]_-{#1}}
\newcommand{\fR}[1]{\ar@{->}[rr]^-{#1}}
\newcommand{\fRr}[1]{\ar@{->}[rrr]^-{#1}}
\newcommand{\fD}[1]{\ar@{->}[dd]_-{#1}}
\newcommand{\fU}[1]{\ar@{->}[uu]^-{#1}}
\newcommand{\f}[2]{\ar@{->}[#1]|{#2}}
\newcommand{\ff}[2]{\ar@2{->}[#1]|{#2}}
\newcommand{\frr}[1]{\ar@{->}[rrrr]^-{#1}}

\newcommand{\fl}[1]{\ar@{->}[l]_-{#1}}
\newcommand{\fr}[1]{\ar@{->}[r]^-{#1}}
\newcommand{\fd}[1]{\ar@{->}[d]_-{#1}}
\newcommand{\fu}[1]{\ar@{->}[u]^-{#1}}

\renewcommand{\top}{{\mathbf{Top}}}

\newcommand{\iso}{\cong}

\DeclareMathOperator{\ev}{ev}

\newcommand{\vI}{\vec{I}}
\renewcommand{\leq}{\leqslant}
\renewcommand{\geq}{\geqslant}

\newcommand{\ptop}[1]{{\brm{{#1}dTop}}}
\newcommand{\tptop}[1]{{\brm{{#1}dTOP}}}

\def\cartesien{%
	\ar@{-}[]+R+<6pt,-2pt>;[]+RD+<6pt,-6pt>%
	\ar@{-}[]+D+<2pt,-6pt>;[]+RD+<6pt,-6pt>%
}
\def\cocartesien{%
	\ar@{-}[]+L+<-6pt,+2pt>;[]+LU+<-6pt,+6pt>%
	\ar@{-}[]+U+<-2pt,+6pt>;[]+LU+<-6pt,+6pt>%
}
\def\hocartesien{%
	\ar@{-}[]+R+<6pt,-2pt>;[]+RD+<6pt,-6pt>_{h}%
	\ar@{-}[]+D+<2pt,-6pt>;[]+RD+<6pt,-6pt>%
}
\def\hococartesien{%
	\ar@{-}[]+L+<-6pt,+2pt>;[]+LU+<-6pt,+6pt>_{h}%
	\ar@{-}[]+U+<-2pt,+6pt>;[]+LU+<-6pt,+6pt>%
}

\newcommand{\brm}[1]{{\rm{\mathbf{#1}}}}

\newcommand{\dtop}{{\brm{Flow}}}
\newcommand{\predtop}{{\mathcal{G}ph(\top)}}
\newcommand{\predtopv}{{\mathcal{G}ph(\mathcal{V})}}

\newcommand{\set}{{\brm{Set}}}

\newcommand{\tdtop}{{\brm{FLOW}}}
\newcommand{\ttop}{{\brm{TOP}}}
\newcommand{\mtop}{{\brm{MSpc}}}

\newcommand{\glob}{{\mathrm{Glob}}}

\DeclareMathOperator{\id}{Id}

\DeclareMathOperator{\Mor}{Mor}

\newcommand{\liminj}{\varinjlim}
\newcommand{\limproj}{\varprojlim}

\renewcommand{\P}{\mathbb{P}}

\makeatletter
\def\varholim@#1#2{%
	\vtop{\m@th\ialign{##\cr
			\hfil$#1\operator@font holim$\hfil\cr
			\noalign{\nointerlineskip\kern1.5\ex@}#2\cr
			\noalign{\nointerlineskip\kern-\ex@}\cr}}%
}
\def\holimproj{%
	\mathop{\mathpalette\varholim@{\leftarrowfill@\textstyle}}\nmlimits@
}
\def\holiminj{%
	\mathop{\mathpalette\varholim@{\rightarrowfill@\textstyle}}\nmlimits@
}
\makeatother

\DeclareMathOperator{\cell}{{\brm{cell}}}
\DeclareMathOperator{\cof}{{\brm{cof}}}
\DeclareMathOperator{\inj}{{\brm{inj}}}

\DeclareMathOperator{\pf}{{Gph}}
\DeclareMathOperator{\cocyl}{{Path}}
\newcommand{\adj}[5]{\xymatrix@C=#5em{{#1}\ar@/^0.8em/[r]^-{#2} \ar@{}[r]|-{\perp} & \ar@/^0.8em/[l]^-{#3} {#4}}}

\begin{document}

\begin{abstract} 
We construct a q-model structure, a h-model structure and a m-model structure on multipointed $d$-spaces and on flows. The two q-model structures are combinatorial and left determined and they coincide with the combinatorial model structures already known on these categories. The four other model structures (the two m-model structures and the two h-model structures) are accessible. We give an example of multipointed $d$-space and of flow which are not cofibrant in any of the model structures. We explain why the m-model structures, Quillen equivalent to the q-model structure of the same category, are better behaved than the q-model structures.
\end{abstract}

\maketitle

\tableofcontents

\section{Introduction}

\subsection*{Presentation}
This paper belongs to our series of papers which aims at comparing the model category $\dtop$ of \textit{flows} introduced in \cite{model3} (with some updated proofs in \cite{leftdetflow} using Isaev's work \cite{Isaev}) and the model category $\ptop{\mathcal{G}}$ of \textit{multipointed $d$-spaces} introduced in \cite{mdtop}. Roughly speaking, the former is a version of the latter without underlying topological space. And the latter is a variant of Grandis' notion of $d$-space \cite{mg}. They are topological models introduced to study concurrent processes from the point of view of homotopy theory. Even if these model categories do not yet contain enough weak equivalences (their homotopical localizations with respect to the refinement of observation remain to be understood: see the digression section in \cite{leftdetflow}), the model category of flows enabled us anyway to understand homological theories detecting the non-deterministic branching and merging areas of execution paths in the framework of flows \cite{3eme} \cite{exbranch}. These homology theories are interesting because they are invariant by the refinement of observation. 

Using the notion of \textit{topological graph} (see Definition~\ref{topgraph}) and the Garner Hess K\k{e}dziorek Riehl Shipley theorem \cite{HKRS17} \cite{GKR18} about accessible right-induced model structures, we introduce a categorical construction which takes as input an accessible model structure on the category $\top$ of $\Delta$-generated spaces satisfying some mild conditions (the ones of Proposition~\ref{preparation-util}) and which gives as output an accessible model structure on multipointed $d$-spaces and on flows. These mild conditions are satisfied in particular by~\footnote{We use the terminology of \cite{ParamHomTtheory}.} the \textit{q-model structure} (the Quillen model structure) of $\top$, the \textit{h-model structure} (also called the Cole-Str\"om model structure) of $\top$ and the \textit{m-model structure} (which is the mixing of the two preceding model structures in the sense of \cite[Theorem~2.1]{mixed-cole}). The latter is characterized as the unique model structure on $\top$ such that the weak equivalences are the weak homotopy equivalences and the fibrations the h-fibrations. We obtain the following results: 
\begin{itemize}[leftmargin=*]
\item a q-model structure, a h-model structure and a m-model structure on multipointed $d$-spaces and on flows \underline{in one step} (!)
\item the identity functor induces a Quillen equivalence between the q-model structure and the m-model structure on multipointed $d$-spaces (on flows resp.)
\item the two q-model structures are combinatorial and left determined and they coincide with that of \cite{mdtop} and of \cite{model3} \cite{leftdetflow} respectively
\item the four other model structures (the two m-model structures and the two h-model structures) are accessible
\item all objects are fibrant in these six model structures
\item there are the implications $\hbox{q-cofibrant} \Rightarrow \hbox{m-cofibrant} \Rightarrow \hbox{h-cofibrant}$ for multipointed $d$-spaces and flows
\item there exist multipointed $d$-spaces and flows which are not q-cofibrant, not h-cofibrant and not m-cofibrant. 
\end{itemize}

The two h-model structures and the two m-model structures are new. They are conjecturally not combinatorial. Even if all topological spaces are h-cofibrant, it is not true that all multipointed $d$-spaces and all flows are h-cofibrant as well. Intuitively, the h-cofibrant objects correspond to objects without algebraic relations in their spaces of execution paths. A rigorous characterization of the h-cofibrant multipointed $d$-spaces and h-cofibrant flows still remains to be find out. 

The main interest of this categorical construction lies in the two m-model structures. They are better behaved than the q-model structures for the following reasons. Unlike the space of execution paths functor $\P:\dtop \to \top$ which preserves q-cofibrancy, it is not true that the space of execution paths functor $\P^{\mathcal{G}}:\ptop{\mathcal{G}}\to \top$ does as well: see Section~\ref{pathspace}. However we have the following result which can be considered as an application of the results of this paper: 

\begin{thmN} (Theorem~\ref{goodbehave} and Theorem~\ref{betterbehave})
The space of execution paths functors $\P^{\mathcal{G}}:\ptop{\mathcal{G}}\to \top$ and $\P:\dtop \to \top$ preserve m-cofibrancy. 
\end{thmN}

We want to end the introduction with a remark about the notion of multipointed $d$-space. It is easy to prove that all theorems of this paper involving multipointed $d$-spaces, except Proposition~\ref{homocat} coming from \cite{model2} and Theorem~\ref{goodbehave}, are still true by replacing the topological group $\mathcal{G}$ of nondecreasing homeomorphisms of the segment $[0,1]$ by the topological monoid $\mathcal{M}$ of nondecreasing continuous maps from the segment $[0,1]$ to itself preserving the extremities. However, we do not know whether Proposition~\ref{homocat} and Theorem~\ref{goodbehave} hold with this new definition of multipointed $d$-space. Indeed, the results of \cite{model2}, in particular Proposition~\ref{homocat} used in the proof of Theorem~\ref{goodbehave}, use the fact that all elements of $\mathcal{G}$ are invertible and we are unable to remove completely this hypothesis by now from the proofs of \cite{model2}.

\subsection*{Outline of the paper}

\begin{itemize}[leftmargin=*]
\item Section~\ref{reminder} collects some basic facts about accessible model categories. It is expounded the theorem we are going to use to right-induce accessible model structures (Theorem~\ref{GHKRS}). 
\item Section~\ref{sorite-bifib} proves two technical elementary facts about Grothendieck bifibrations that will be used in the sequel: a first one which is a toolkit to easily prove that a functor is a bifibration (Proposition~\ref{toolbifib}), and a second one about the accessibility of two functors arising from an accessible bifibration (Proposition~\ref{munuaccessible}).
\item Section~\ref{topspace} gathers some information about $\Delta$-generated spaces and their three standard model structures. In particular, Proposition~\ref{preparation-util} makes explicit and establishes that these three model structures satisfy the mild conditions which are used in our construction.
\item Section~\ref{enrichedgraph} explains how to construct an accessible model structure on $\mathcal{V}$-graphs from any accessible model category $\mathcal{V}$ (Theorem~\ref{any}), with an immediate application when $\mathcal{V}$ is the category of $\Delta$-generated spaces (Corollary~\ref{qhm-graph}).
\item Section~\ref{homotopy-dspace} applies the constructions of Section~\ref{enrichedgraph} to right-induce on the category of multipointed $d$-spaces the three model structures (Theorem~\ref{three}). It is also proved that there exist multipointed $d$-spaces which are not h-cofibrant, not q-cofibrant and not m-cofibrant (Proposition~\ref{mdtop-notcofibrant}).
\item Section~\ref{homotopy-flow} applies the same constructions to right-induce on the category of flows the three model structures (Theorem~\ref{three2}). It is also proved that there exist flows which are not h-cofibrant, not q-cofibrant and not m-cofibrant (Proposition~\ref{flow-notcofibrant}).
\item Section~\ref{pathspace} explains why the m-model structures are better behaved than the q-model structures (Theorem~\ref{goodbehave} and Theorem~\ref{betterbehave}).
\end{itemize}

\subsection*{Notations}

\begin{itemize}[leftmargin=*]
\item $X:=Y$ means that $Y$ is the definition of $X$.
\item All categories are locally small (except the category of all locally small categories).
\item $\K$ always denotes a locally presentable category.
\item $\set$ is the category of sets.
\item $\top$ is the category of $\Delta$-generated spaces.
\item $\mathcal{G}$ is the topological group of nondecreasing homeomorphisms of $[0,1]$.
%\item $\mathcal{M}$ is the topological monoid of nondecreasing continuous maps from $[0,1]$ to itself preserving the extremities.
\item $\mathbb{R}$ is the topological space of real numbers.
\item $\K(X,Y)$ is the set of maps in a category $\K$.
\item $\Mor(\K)$ is the category of morphisms of $\K$ with the commutative squares for the morphisms.
\item $A \sqcup B$ is the binary coproduct, $A \p B$ is the binary product. 
\item $\limproj$ is the limit, $\liminj$ is the colimit.
\item $\varnothing$ is the initial object.
\item $\mathbf{1}$ is the final object.
\item $\id_X$ is the identity of $X$.
\item $g.f$ is the composite of two maps $f:A\to B$ and $g:B\to C$; the composite of two functors is denoted in the same way. 
\item $f\boxslash g$ means that $f$ satisfies the \textit{left lifting property} (LLP) with respect to $g$, or equivalently that $g$ satisfies the \textit{right lifting property} (RLP) with respect to $f$.
\item $\inj(\C) = \{g \in \K, \forall f \in \C, f\boxslash g\}$.
\item $\cof(\C)=\{f\mid \forall g\in \inj(\C), f\boxslash g\}$.
\item $\cell(\C)$ is the class of transfinite compositions of pushouts of elements of $\C$.
\item A \textit{cellular} object $X$ of a combinatorial model category is an object such that the canonical map $\varnothing\to X$ belongs to $\cell(I)$ where $I$ is the set of generating cofibrations.
\item A \textit{model structure} $(\C,\W,\F)$ means that the class of cofibrations is $\C$, that the class of weak equivalences is $\W$ and that the class of fibrations is $\F$ in this order. A \textit{model category} is a category equipped with a \textit{model structure}.
\item $\ell,\ell_i,\ell',\ell'_i$ always denote nonzero positive real numbers. 
\item The notation $[0,\ell_1]\iso^+ [0,\ell_2]$ means a nondecreasing homeomorphism from $[0,\ell_1]$ to $[0,\ell_2]$. It takes $0$ to $0$ and $\ell_1$ to $\ell_2$. The group for the composition of maps of nondecreasing homeomorphisms from $[0,1]$ to itself is denoted by $\mathcal{G}$, i.e. $\mathcal{G}=\{[0,1]\iso^+[0,1]\}$.
%\item The notation $[0,\ell_1]\to^+ [0,\ell_2]$ means a nondecreasing continuous from $[0,\ell_1]$ to $[0,\ell_2]$ preserving the extremities. It takes $0$ to $0$ and $\ell_1$ to $\ell_2$. The monoid for the composition of maps of nondecreasing continuous maps from $[0,1]$ to itself preserving the extremities is denoted by $\mathcal{M}$, i.e. $\mathcal{M}=\{[0,1]\to^+[0,1]\}$.
%\item The letter $\mathcal{P}$ denotes $\mathcal{G}$ or $\mathcal{M}$. 
%
%
%
\end{itemize}

\subsection*{Acknowledgments}

I thank the MathOverflow community. Many contributions drew my attention to the recent breakthroughs in the theory of accessible model categories. I also thank the contributors of the nLab website. I thank the referee for the helpful comments.

\section{Accessible model category}
\label{reminder}

We refer to \cite{TheBook} for locally presentable categories, to \cite{MR2506258} for combinatorial model categories.  We refer to \cite{MR99h:55031} and to \cite{ref_model2} for more general model categories. 

A weak factorization system $(\mathcal{L},\mathcal{R})$ of a locally presentable category $\K$ is \textit{accessible} if there is a functorial factorization 
\[\xymatrix@1{(A\stackrel{f}\longrightarrow B) \ar@{|->}[r] & (A\stackrel{Lf}\longrightarrow Ef\stackrel{Rf}\longrightarrow B)}\]
with $Lf\in \mathcal{L}$, $Rf\in \mathcal{R}$ such that the functor $E:\Mor(\K)\to \K$ is accessible \cite[Definition~2.4]{GKR18}. Since colimits are calculated pointwise in $\Mor(\K)$, a weak factorization system is accessible if and only if the functors $L:\Mor(\K)\to \Mor(\K)$ and $R:\Mor(\K)\to \Mor(\K)$ are accessible. By \cite[Theorem~4.3]{MR3638359}, a weak factorization system is accessible if and only if it is small in Garner's sense. In particular, every \textit{small} weak factorization system (i.e. of the form $(\cof(I),\inj(I))$ for a set $I$) is accessible. A model structure $(\C,\W,\F)$ on a locally presentable category is \textit{accessible} if the two weak factorization systems $(\C,\W\cap\F)$ and $(\C\cap\W,\F)$ are accessible. Every combinatorial model category is therefore an accessible model category. This inclusion is strict: by \cite[Remark~4.7]{Raptis-Strom}, the h-model structure of $\top$ (see Section~\ref{topspace}) is not combinatorial. But it is accessible by Proposition~\ref{preparation-util}. Moreover, there exist model categories which are not accessible. For example, the model category of maps of spaces of \cite{ModelCategoryMapsOfSpaces} is conjecturally not accessible (remark due to Boris Chorny \cite{Ex1}) but no proof has been given yet. There is another example using the negation of Vop\v{e}nka's principle given by Mike Shulman \cite{Ex2}: by \cite[Example~6.12]{TheBook}, the locally presentable category $\mathbf{Gra}$ of graphs has a reflective subcategory that is not accessible if we assume the negation of Vop\v{e}nka's principle, and by \cite[Proposition~3.5]{DiscreteModelStructure}, this reflector is the fibrant replacement functor of a model structure on $\mathbf{Gra}$. 

The following theorem is the particular case of a general theorem due to Garner Hess K\k{e}dziorek Riehl and Shipley about accessible right-induced model structures (note that the Quillen Path Object argument dates back to \cite{MR36:6480}). 

\bth \label{GHKRS} (Garner-Hess-K\k{e}dziorek-Riehl-Shipley) Let $\mathcal{M}$ and $\mathcal{N}$ be two locally presentable categories. Let $(\C,\W,\F)$ be an accessible model structure of $\mathcal{M}$ such that all objects are fibrant. Consider a categorical adjunction \[\adj{\mathcal{M}}{L}{U}{\mathcal{N}}{5}.\] Suppose that there exists a functorial factorization of the diagonal of $\mathcal{N}$
\[\xymatrix@C=3em{X \fr{\tau}& \cocyl(X) \fr{\pi}& X\p X}\]
such that $U(\tau)$ is a weak equivalence of $\mathcal{M}$ and such that $U(\pi)$ is a fibration of $\mathcal{M}$ for all objects $X$ of $\mathcal{N}$. Then there exists a unique model structure on $\mathcal{N}$ such that the class of fibrations is $U^{-1}(\F)$ and such that the class of weak equivalences is $U^{-1}(\W)$. Moreover, this model structure is accessible and all its objects are fibrant. 
\eth

\bpf[Sketch of proof] By the dual of \cite[Theorem~2.2.1]{HKRS17} which is also stated in \cite[Theorem~6.2]{MoserLyne}, the hypotheses of the theorem imply that the Quillen Path Object argument holds. The latter implies the acyclicity condition for right-induced model structures, and therefore the existence of the right-induced model structure (see also \cite{GKR18}). Since a model structure is characterized by its class of weak equivalences and its class of fibrations, we deduce the uniqueness.
\epf

\section{Accessible Grothendieck bifibration}
\label{sorite-bifib}

Let $p:\mathcal{E}\to \mathcal{B}$ be a functor between locally small categories. The \textit{fibre of $p$ over $X$}, denoted by $\mathcal{E}_X$, consists of the subcategory of $\mathcal{E}$ generated by the \textit{vertical maps} $f$, i.e. the maps $f$ such that $p(f)=\id_X$. We refer to \cite[Chapter~1 and Chapter~9]{JacobBook} and \cite[Chapter~8]{Borceux2} for \textit{(Grothendieck) bifibrations} (also called \textit{bifibred categories}) and  for \textit{(Grothendieck) fibrations} (also called \textit{fibred categories}, the term fibration being quite confusing because it is used in a completely different sense in this paper).

The following proposition is a toolkit to minimize the work required to prove that a functor is a bifibration: 

\bp \label{toolbifib}
Let $p:\mathcal{E}\to \mathcal{B}$ be a functor between locally small categories. Suppose that for every map $u:A\to B$ of $\mathcal{B}$, there exists an adjunction $u_!:\mathcal{E}_A \dashv \mathcal{E}_B:u^*$ such that: 
\begin{enumerate}[leftmargin=*]
\item For all objects $X$ of $\mathcal{E}$, there exists a natural map $u^*X\to X$ such that every map $f:X\to Y$ of $\mathcal{E}$ with $p(f)=u$ factors uniquely as a composite \[X\longrightarrow u^*Y\longrightarrow Y\] with the left-hand map vertical.
\item The natural map $u^*v^*X\to (v.u)^*X$ is an isomorphism for all $X$.
\end{enumerate}
Then $p:\mathcal{E}\to \mathcal{B}$ is a bifibration.
\ep

\bpf
In the language of \cite{JacobBook}, the first condition means that the map $u^*X\to X$ is weakly cartesian and the second condition implies that compositions of weakly cartesian maps are weakly cartesian. By \cite[Exercice~1.1.6]{JacobBook}, the functor $p:\mathcal{E}\to \mathcal{B}$ is a fibred category. By \cite[Lemma~9.1.2]{JacobBook}, the existence of the adjunctions implies that the functor $p:\mathcal{E}\to \mathcal{B}$ is a bifibration. 
\epf

Let $p:\mathcal{E}\to \mathcal{B}$ be a bifibration between locally small categories. Consider the commutative square of solid arrows of $\mathcal{E}$
	\[
	\xymatrix@C=3em@R=3em
	{
		X \ar@/^20pt/@{->}[rr]^-{f}\fr{}\fd{g} & \mu(f):=p(f)^*Y \ar@{-->}[d]_-{\exists !}^-{\mu_{(g,h)}}\fr{} & Y \ar@{->}[d]^-{h}\\
		X'\ar@/_20pt/@{->}[rr]_-{f'} \fr{} & \mu(f'):=p(f')^*Y' \fr{} & Y'.
	}
	\]
Note that the diagram above is misleading: the maps $g$ and $h$ are not vertical. On the contrary, the two maps $X\to \mu(f)$ and $X'\to \mu(f)$ are vertical. Since $\mu(f')\to Y'$ is cartesian, there exists a unique map $\mu_{(g,h)}:\mu(f) \to \mu(f')$ such that $p(\mu_{(g,h)}) = p(g)$ making the right-hand square commutative. Since the composites $X\to \mu(f)\to \mu(f')$ and $X\to X'\to \mu(f')$ have the same image $p(g)$ by $p$ and since they yield two factorizations of $h.f=f'.g$ and since $\mu(f')\to Y'$ is cartesian, the left-hand square is commutative as well. For dual reasons, there exists a unique map $\nu_{(g,h)}:\nu(f)\to \nu(f')$ such that $p(\nu_{(g,h)}) = p(h)$ making the following diagram of solid arrows of $\mathcal{E}$
	\[
	\xymatrix@C=3em@R=3em
	{
		X \ar@/^20pt/@{->}[rr]^-{f}\fr{}\fd{g} & \nu(f):=p(f)_!X \ar@{-->}[d]_-{\exists !}^-{\nu_{(g,h)}}\fr{} & Y \ar@{->}[d]^-{h}\\
		X'\ar@/_20pt/@{->}[rr]_-{f'} \fr{} & \nu(f'):=p(f')_!X' \fr{} & Y'.
	}
	\]
commutative. By the usual uniqueness argument, we obtain two well-defined functors $\mu:\Mor(\mathcal{E})\to \mathcal{E}$ and $\nu:\Mor(\mathcal{E})\to \mathcal{E}$.

\bp \label{munuaccessible} Let $p:\mathcal{E}\to \mathcal{B}$ be a bifibration between locally presentable categories such that $p$ is accessible. Then the functors $\mu:\Mor(\mathcal{E})\to \mathcal{E}$ and $\nu:\Mor(\mathcal{E})\to \mathcal{E}$ defined above are accessible.
\ep

\bpf Suppose that $p:\mathcal{E}\to \mathcal{B}$ is $\lambda$-accessible. Let $(f_i:X_i\to Y_i)_{i\in I}$ be a $\lambda$-filtered diagram of $\Mor(\mathcal{E})$. By passing to the colimit, we obtain the factorization of $\liminj f_i$ 
\[
\liminj X_i \longrightarrow \liminj \mu(f_i) \longrightarrow \liminj Y_i.
\]
There are the isomorphisms \[p(\liminj X_i)\iso \liminj p(X_i) = \liminj p(\mu(f_i)) \iso p(\liminj \mu(f_i)),\] the first isomorphism since $p$ is $\lambda$-accessible, the equality since each $X_i\to \mu(f_i)$ is vertical, and the last isomorphism since $p$ is $\lambda$-accessible. Let $u:p(\liminj X_i)\to p(\liminj \mu(f_i))$ be this isomorphism. Then we have the isomorphism \[u^*(\liminj \mu(f_i)) \iso \liminj \mu(f_i).\] We obtain the factorization of $\liminj f_i$ 
\[
\liminj X_i \longrightarrow u^*(\liminj \mu(f_i)) \longrightarrow \liminj Y_i.
\]
Since the left-hand map is vertical, we obtain the equality \[\mu(\liminj f_i) = u^*(\liminj \mu(f_i)).\] We have proved that $\mu$ is accessible. In the same way, by passing to the colimit, there is the factorization of $\liminj f_i$ 
\[
\liminj X_i \longrightarrow \liminj \nu(f_i) \longrightarrow \liminj Y_i. 
\]
There are the isomorphisms \[p(\liminj \nu(f_i)) \iso \liminj p(\nu(f_i)) = \liminj p(Y_i) \iso p(\liminj Y_i),\] the first isomorphism since $p$ is $\lambda$-accessible, the equality since each $\nu(f_i) \to Y_i$ is vertical and the last isomorphism since $p$ is $\lambda$-accessible. Let $v:p(\liminj \nu(f_i))\to p(\liminj Y_i)$ be this isomorphism. Then we have the isomorphism \[v_!(\liminj \nu(f_i)) \iso \liminj \nu(f_i).\] We obtain the factorization of $\liminj f_i$ 
\[
\liminj X_i \longrightarrow v_!(\liminj \nu(f_i)) \longrightarrow \liminj Y_i. 
\]
Since the right-hand map is vertical, we obtain the equality \[\nu(\liminj f_i) = v_!(\liminj \nu(f_i)).\] We have proved that $\nu$ is accessible.
\epf

\section{Delta-generated space}
\label{topspace}

We refer to \cite[Chapter~VI]{topologicalcat} or \cite[Chapter~7]{Borceux2} for the notion of topological functor. The category $\top$ denotes the category of \textit{$\Delta$-generated spaces}, i.e. the colimits of simplices. Let $\Delta^n=\{(t_0,\dots,t_n)\in [0,1]^n\mid t_0+\dots +t_n=1\}$ be the topological $n$-simplex equipped with its standard topology. Then $\top$ is the final closure of the set of topological spaces $\{\Delta^n\mid n\geq 0\}$. For a tutorial about these topological spaces, see for example \cite[Section~2]{mdtop}. The category $\top$ is locally presentable by \cite[Corollary~3.7]{FR}, cartesian closed and it contains all CW-complexes. The internal hom functor is denoted by $\ttop(-,-)$. We denote by $\omega:\mathcal{T\!O\!P}\to \set$ the underlying set functor where $\mathcal{T\!O\!P}$ is the category of general topological spaces. It is fibre-small and topological. The restriction functor $\omega:\top\subset \mathcal{T\!O\!P}\to \set$ is fibre-small and topological as well. The category $\top$ is a full coreflective subcategory of the category $\mathcal{T\!O\!P}$ of general topological spaces. Let $k:\mathcal{T\!O\!P}\to\top$ be the kelleyfication functor, i.e. the right adjoint. The category $\top$ is finally closed in $\mathcal{T\!O\!P}$, which means that the final topology and the $\omega$-final structure coincides. On the contrary, the $\omega$-initial structure in $\top$ is obtained by taking the kelleyfication of the initial topology in $\mathcal{T\!O\!P}$. If $A$ is a subset of a space $X$ of $\top$, the initial structure in $\top$ of the inclusion $A\subset \omega X$ is the kelleyfication of the relative topology with respect to the inclusion.

\begin{rem}
It is important to keep in mind for the sequel that the kelleyfication functor does not change the underlying set. In particular, it does not identify points. It only adds open sets to the topology.
\end{rem}

\begin{nota}
Let $n\geq 1$. Denote by $\mathbf{D}^n = \{b\in \mathbb{R}^n, |b| \leq 1\}$ the $n$-dimensional disk, and by $\mathbf{S}^{n-1} = \{b\in \mathbb{R}^n, |b| = 1\}$ the $(n-1)$-dimensional sphere. By convention, let $\mathbf{D}^{0}=\{0\}$ and $\mathbf{S}^{-1}=\varnothing$. 
\end{nota}

The category $\top$ can be equipped at least with three model structures (we use the notations of \cite{ParamHomTtheory}): 
\begin{itemize}[leftmargin=*]
\item The \textit{q-model structure} $(\C_q,\W_q,\F_q)$ \cite[Section~2.4]{MR99h:55031}: the cofibrations, called \textit{q-cofibrations}, are the retracts of the transfinite compositions of the inclusions $\mathbf{S}^{n-1}\subset \mathbf{D}^n$ for $n\geq 0$, the weak equivalences are the weak homotopy equivalences and the fibrations, called \textit{q-fibrations} are the maps satisfying the RLP with respect to the inclusions $\mathbf{D}^{n}\subset \mathbf{D}^{n+1}$ for $n\geq 0$, or equivalently with respect to the inclusions $\mathbf{D}^{n}\p \{0\}\subset \mathbf{D}^{n}\p [0,1]$ for $n\geq 0$; this model structure is combinatorial. A very simple way to obtain this model structure is to use \cite{Isaev}. Its existence dates back to \cite{MR36:6480}.
\item The \textit{h-model structure} $(\C_{\overline{h}},\W_h,\F_h)$: the fibrations, called the \textit{h-fibrations}, are the maps satisfying the RLP with respect to the inclusions $X\p \{0\}\subset X\p [0,1]$ for all topological spaces $X$, and the weak equivalences are the homotopy equivalences;  we have $\C_q\subset \C_{\overline{h}}$ because $\W_h\subset \W_q$ and $\F_h\subset \F_q$. A modern exposition is given in \cite[Corollary~5.23]{Barthel-Riel} but its construction dates back to \cite{strom3}. All topological spaces are h-cofibrant.
\item The \textit{m-model structure} $(\C_m,\W_m,\F_m)=(\C_m,\W_q,\F_h)$: the fibrations are the \textit{h-fibra\-tions}, and the weak equivalences are the weak homotopy equivalences; we have $\C_q\subset \C_m$ because $\W_m\cap \F_m=\W_q\cap \F_h \subset \W_q\cap \F_q$. Its existence is a consequence of \cite[Theorem~2.1]{mixed-cole}. By \cite[Corollary~3.7]{mixed-cole}, a topological space is m-cofibrant if and only if it is homotopy equivalent to a q-cofibrant space.
\end{itemize}

\begin{figure}
\[\xymatrix{X\iso \{0\}\p X \fr{f}\fd{\subset} & Y \fd{} &\Pi Y \cartesien\fr{}\fd{}& \ttop(\mathbb{R}^+,Y)\p \mathbb{R}^+ \ar@{->}[d]^-{\hbox{\tiny shift}}   & \Gamma f\fd{} \cartesien\fr{}& \Pi Y\fd{p_0}\\
[0,1]\p X \fr{} & \cocartesien Mf & Y \ar@{->}[r]_-{{\hbox{\tiny const}}} & \ttop(\mathbb{R}^+,Y)  & X \fr{f} & Y.
}\]
\caption{Mapping cylinder $Mf$ and Moore path space $\Gamma f$ with $\mathbb{R}^+=[0,+\infty[$}
\label{mapping}
\end{figure}

\bp \label{preparation-util} The three model structures of $\top$ $(\C_q,\W_q,\F_q)$, $(\C_{\overline{h}},\W_h,\F_h)$ and $(\C_m,\W_m,\F_m)$ satisfy the following properties: 
\begin{enumerate}
\item They are accessible.
\item All spaces are fibrant.
\item All homotopy equivalences are weak equivalences.
\item All q-cofibrations are cofibrations.
\item For all topological spaces $X$ of $\top$, the map \[\pi:\ttop([0,1],X)\to \ttop(\{0,1\},X)\] induced by the inclusion $\{0,1\}\subset [0,1]$ is a fibration.
\end{enumerate}
\ep

\bpf $(1)$ The model structure $(\C_q,\W_q,\F_q)$ is accessible because it is combinatorial. The model structure $(\C_m,\W_m,\F_m)$ is accessible by \cite[Corollary~4.4]{dgrtop}. Figure~\ref{mapping} recalls the definition of the Moore paths space $\Pi Y$ of $Y$ of \cite[Section~3.1]{Barthel-Riel} which actually dates back to \cite{MR0370579}. The bottom map $Y\to \ttop(\mathbb{R}^+,Y)$ is the constant path map. The shift map $\ttop(\mathbb{R}^+,Y)\p \mathbb{R}^+\to \ttop(\mathbb{R}^+,Y)$ takes the pair $(\gamma,t)$ to the path $u\mapsto \gamma(t+u)$. By definition $p_0(\gamma)=\gamma(0)$. Since $\top$ is locally presentable and cartesian closed, it is easily seen that the Moore path functor $\Gamma:f\mapsto \Gamma f$ of Figure~\ref{mapping} is accessible. It provides a functorial factorization for $(\C_{\overline{h}}\cap\W_h,\F_h)$ by \cite[Corollary~3.12]{Barthel-Riel}. The functorial factorization $(\C_{\overline{h}},\W_h\cap\F_h)$ is given first by taking a map $f:X\to Y$ of $\top$ to the composite map $X\to M f \to Y$ and then by using on the right-hand map the functorial factorization of $(\C_{\overline{h}}\cap\W_h, \F_h)$ using the Moore path functor (see \cite[Proposition~3.2]{MR1967263}). This proves that the model structure $(\C_{\overline{h}},\W_h,\F_h)$ is accessible. $(2)$ and $(3)$ are well-known. $(4)$ is recalled above. $(5)$ deserves a short bibliographical justification. The inclusion $\{0,1\}\subset [0,1]$ is a cofibration in the three model structures. And the canonical map $X\to \mathbf{1}$ is a fibration in the three model structures as well. It suffices to use \cite[Lemma~4.2.2(3)]{MR99h:55031} and the fact that the three model structures are monoidal for the binary product: for the q-model structure of $\top$, see e.g. \cite[Proposition~4.2.11]{MR99h:55031}; for the h-model structure of $\top$, see e.g. \cite[Corollary~2.10]{MR1967263} for a general treatment in the setting of enriched categories; for the m-model structure, see e.g. \cite[Proposition~6.6]{mixed-cole}.
\epf

\section{Topological graph}
\label{enrichedgraph}

In this section, $\mathcal{V}$ denotes a locally presentable category. It is supposed to be equipped with an accessible model structure $(\C,\W,\F)$. We recall the enriched version of the usual notion of graph and of morphism between them \cite[Definition~5.1.1]{Borceux1}. This notion appears for example in \cite[Definition 2.1.1]{EnrichedGraph} and in \cite[Section~3]{VCatLocPre}. We adapt the notations to our context.

\bd A {\rm $\mathcal{V}$-graph} $X$ consists of a pair \[(X^0,(\P_{\alpha,\beta}X)_{(\alpha,\beta)\in X^0\p X^0})\] such that $X^0$ is a set and such that each $\P_{\alpha,\beta}X$ is an object of $\mathcal{V}$. A map of $\mathcal{V}$-graphs $f:X\to Y$ consists of a set map $f^0:X^0\to Y^0$ (called the {\rm underlying set map}) together with a map $\P_{\alpha,\beta}X\to \P_{f^0(\alpha),f^0(\beta)}Y$ of $\mathcal{V}$ for all $(\alpha,\beta)\in X^0\p X^0$. The composition is defined in an obvious way. The corresponding category is denoted by $\predtopv$. 
\ed

\begin{nota}
We will denote $\P_{f^0(\alpha),f^0(\beta)}Y$ by $\P_{f(\alpha),f(\beta)}Y$ in order not to overload the notations.
\end{nota}

\bp \label{fibred} The forgetful functor $X\mapsto X^0$ from $\predtopv$ to $\set$ is a bifibration. \ep

\bpf Let $f:X\to Y$ be a map of $\mathcal{V}$-graphs. Let 
\[\P_{\alpha,\beta} (f^0)^*Y := \P_{f(\alpha),f(\beta)}Y\]
for all $(\alpha,\beta)\in X^0\p X^0$. We obtain a well-defined $\mathcal{V}$-graph $(f^0)^*Y$. Then by definition of a map of $\mathcal{V}$-graphs, every map $f:X\to Y$ factors uniquely as a composite \[X\longrightarrow (f^0)^*Y \longrightarrow Y\] such that the left-hand map is vertical. Thus the map $(f^0)^*Y \to Y$ is weakly cartesian. The fact that $(g^0.f^0)^* = (f^0)^*.(g^0)^*$ for two composable maps $f$ and $g$ is obvious. Let 
\[
\P_{\gamma,\delta}(f^0)_!X  = \bigsqcup_{{\tiny\begin{array}{c}
(\alpha,\beta)\in X^0\p X^0 \\f(\alpha)=\gamma,f(\beta)=\delta
\end{array}}} \P_{\alpha,\beta}X.
\]
for all $(\gamma,\delta)\in Y^0\p Y^0$. We obtain a well-defined $\mathcal{V}$-graph $(f^0)_!X$. We have the natural bijections of sets 
\begin{align*}
\predtopv_{X^0}(X,(f^0)^*Y)  &\iso \prod_{(\alpha,\beta)\in X^0\p X^0} \mathcal{V}(\P_{\alpha,\beta}X,\P_{f(\alpha),f(\beta)}Y) \\
& \iso \prod_{(\gamma,\delta)\in Y^0\p Y^0} \prod_{{\tiny\begin{array}{c}
(\alpha,\beta)\in X^0\p X^0 \\f(\alpha)=\gamma,f(\beta)=\delta
\end{array}}} \mathcal{V}(\P_{\alpha,\beta}X,\P_{\gamma,\delta}Y)\\
& \iso \prod_{(\gamma,\delta)\in Y^0\p Y^0} \mathcal{V}(\P_{\gamma,\delta}(f^0)_!X,\P_{\gamma,\delta}Y)\\
&\iso \predtopv_{Y^0}((f^0)_!X,Y),
\end{align*} 
the first and the fourth isomorphisms by definition of a map of $\mathcal{V}$-graphs, the second isomorphism by rearranging the product and the third isomorphism by definition of the $\mathcal{V}$-graph $(f^0)_!X$. The proof is complete thanks to Proposition~\ref{toolbifib}.
\epf

For every set $S$, the fibre of $()^0:\predtopv\to \set$ over $S$ is the functor category $\mathcal{V}^{S\p S}$ which is equipped for the sequel with the only model structure such that the cofibrations (the fibrations, the weak equivalences resp.) are the pointwise ones: it is both the projective and the injective model structure on a functor category over a discrete category. This model structure is obviously accessible. 

\bth \label{any} There exists a unique model structure on $\predtopv$ such that 
\begin{itemize}
\item The weak equivalences are the maps of $\mathcal{V}$-graphs $f:X\to Y$ such that $f^0$ is a bijection and such that the map $X\to (f^0)^*Y$ is a pointwise weak equivalence of $\mathcal{V}^{X^0\p X^0}$, i.e. for all $(\alpha,\beta)\in X^0\p X^0$, the map $\P_{\alpha,\beta}X\to \P_{f(\alpha),f(\beta)}Y$ belongs to $\W$.
\item The fibrations are the maps of $\mathcal{V}$-graphs $f$ such that the map $X\to (f^0)^*Y$ is a pointwise fibration of $\mathcal{V}^{X^0\p X^0}$, i.e. for all $(\alpha,\beta)\in X^0\p X^0$, the map $\P_{\alpha,\beta}X\to \P_{f(\alpha),f(\beta)}Y$ belongs to $\F$.
\item The cofibrations are the maps of $\mathcal{V}$-graphs $f$ such that the map $(f^0)_!X\to Y$ is a pointwise cofibration of $\mathcal{V}^{Y^0\p Y^0}$, i.e. for all $(\gamma,\delta)\in Y^0\p Y^0$, the map $\bigsqcup_{{\tiny\begin{array}{c}
(\alpha,\beta)\in X^0\p X^0 \\f(\alpha)=\gamma,f(\beta)=\delta
\end{array}}} \P_{\alpha,\beta}X \to \P_{\gamma,\delta} Y$ belongs to $\C$. 
\end{itemize}
Moreover, this model structure is accessible. 
\eth

\bpf We want to apply \cite[Theorem~5.1]{Roig} fixed in \cite[Theorem page 23]{Roig-fixed} to the bifibration $()^0:\predtopv\to \set$. We equip the base category $\set$ with the discrete model structure: all maps are cofibrations and fibrations and the weak equivalences are the bijections. For every set map $u:S\to T$, the functor $u^*:\mathcal{V}^{T\p T}\to \mathcal{V}^{S\p S}$ preserves weak equivalences and fibrations since they are pointwise. Therefore, the adjunction $(u_!,u^*)$ is a Quillen adjunction. We have to verify the two hypotheses of \cite[Theorem page 23]{Roig-fixed}: 
\begin{enumerate}[leftmargin=*]
\item if $u:S\to T$ is a weak equivalence of $\set$, then it is a bijection. Therefore the functor $u^*:\mathcal{V}^{T\p T}\to \mathcal{V}^{S\p S}$ reflects weak equivalences since it is an equivalence of categories. 
\item if $u:S\to T$ is a trivial cofibration of $\set$, then it is a bijection, which means that we can suppose that $S=T$. In that case, both $u_!$ and $u^*$ are the identity of $\mathcal{V}^{S\p S}$ and the unit of the adjunction $X\to u^*u_!X$ is an isomorphism, and therefore a weak equivalence of $\mathcal{V}^{S\p S}$.
\end{enumerate}
This proves the existence of the model structure. By \cite[Proposition~4.4]{VCatLocPre}, the category $\predtopv$ is locally presentable~\footnote{This can be proved directly by observing that the fibred category $(-)^0:\predtopv\to\set$ corresponds to an accessible pseudo-functor in the sense of \cite[Definition~5.3.1]{MR1031717} and by applying \cite[Theorem~5.3.4]{MR1031717}.}. Let $f:X\to Y$ be a map of $\mathcal{V}$-graphs. It factors as a composite 
\[
\xymatrix@C=3em
{
X \ar@{^(->}[r]^{\simeq} & Z \ar@{->>}[r] & \mu(f) \fr{} & Y
}
\]
where the factorization trivial cofibration-fibration of the vertical map $X\to \mu(f)$ is carried out in $\mathcal{V}^{X^0\p X^0}$. Since the map $Z\to \mu(f)$ is vertical, we have \[\mu(Z\to Y)=\mu(f)=(f^0)^*Y.\] Thus the composite $Z\to \mu(f)\to Y$ is a fibration of $\predtopv$ by definition of them. We have obtained a factorization trivial cofibration-fibration in $\predtopv$. The functor $(-)^0:\predtopv \to \set$ is colimit preserving since it has a right adjoint: the functor taking a set $S$ to the constant diagram $\Delta_{S\p S}(\mathbf{1})$ over $S\p S$. By Proposition~\ref{munuaccessible}, the endofunctor of $\Mor(\predtopv)$ taking $f:X\to Y$ to $X\to \mu(f)$ is accessible since colimits are calculated pointwise in $\Mor(\predtopv)$. Since the model structure of $\mathcal{V}^{X^0\p X^0}$ is accessible, we deduce that the factorization trivial cofibration-fibration in $\predtopv$ is accessible. The map $f:X\to Y$ factors as well as a composite 
\[
\xymatrix@C=3em
{
X \fr{} & \nu(f) \ar@{^(->}[r] & T \ar@{->>}[r]^{\simeq} & Y
}
\]
where the factorization cofibration-trivial fibration of the vertical map $\nu(f)\to Y$ is carried out in $\mathcal{V}^{Y^0\p Y^0}$. Since the map $\nu(f)\to T$ is vertical, we have \[\nu(X\to T)=\nu(f)=(f^0)_!X.\] Thus the composite $X\to \nu(f)\to T$ is a cofibration of $\predtopv$ by definition of them. We have obtained a factorization cofibration-trivial fibration in $\predtopv$. Since colimits of maps are calculated pointwise, we deduce that the endofunctor of $\Mor(\predtopv)$ taking $f:X\to Y$ to $\nu(f)\to Y$ is accessible by Proposition~\ref{munuaccessible}. Since the model structure of $\mathcal{V}^{Y^0\p Y^0}$ is accessible, we deduce that the factorization cofibration-trivial fibration in $\predtopv$ is accessible. We have proved that the model category $\predtopv$ is an accessible model category.
\epf

\bd \label{topgraph} A {\rm topological graph} is a $\mathcal{V}$-graph with $\mathcal{V}=\top$. The corresponding category is denoted by $\predtop$. \ed

\begin{cor} \label{qhm-graph} Let $(\C,\W,\F)$ be one of the three model structures \[(\C_q,\W_q,\F_q), (\C_{\overline{h}},\W_h,\F_h),(\C_m,\W_m,\F_m)\] of $\top$. Then there exists a unique model structure on $\predtop$ such that: 
\begin{itemize}
\item A map of topological graphs $f:X\to Y$ is a weak equivalence if and only if $f^0:X^0\to Y^0$ is a bijection and for all $(\alpha,\beta)\in X^0\p X^0$, the continuous map $\P_{\alpha,\beta}X\to \P_{f(\alpha),f(\beta)}X$ belongs to $\W$.
\item A map of topological graphs $f:X\to Y$ is a fibration if and only if for all $(\alpha,\beta)\in X^0\p X^0$, the continuous map $\P_{\alpha,\beta}X\to \P_{f(\alpha),f(\beta)}X$ belongs to $\F$.
\end{itemize}
Moreover, this model structure is accessible and all objects are fibrant.
\end{cor}

\bpf
It is a consequence of Theorem~\ref{any} and Proposition~\ref{preparation-util} (1) and (2). 
\epf

\section{Multipointed d-space}
\label{homotopy-dspace}

\bd A {\rm multipointed space} is a pair $(|X|,X^0)$ where
\begin{itemize}
\item $|X|$ is a topological space called the {\rm underlying space} of $X$.
\item $X^0$ is a subset of $|X|$ called the {\rm set of states} of $X$.
\end{itemize}
A morphism of multipointed spaces $f:X=(|X|,X^0) \rightarrow Y=(|Y|,Y^0)$ is a commutative square
\[
\xymatrix{
X^0 \fr{f^0}\fd{} & Y^0 \fd{} \\ 
|X| \fr{|f|} & |Y|.}
\] 
The corresponding category is denoted by $\mtop$.  \ed

For any topological space $U$, two continuous maps $\gamma_1:[0,\ell_1]\to U$ and $\gamma_2:[0,\ell_2]\to U$ with $\ell_1,\ell_2>0$ are \textit{composable} if $\gamma_1(\ell_1) = \gamma_2(0)$. Then one can define the continuous map $\gamma_1*\gamma_2:[\ell_1+\ell_2] \to U$ by 
\[
(\gamma_1 *\gamma_2)(t) = \begin{cases}
\gamma_1(t) & \hbox{if } t\in [0,\ell_1]\\
\gamma_2(t - \ell_1)& \hbox{if }t\in [\ell_1,\ell_1+\ell_2].
\end{cases}
\]
If $\gamma_3:[0,\ell_3]\to U$ is a third continuous map, then there is the (strict) equality \[(\gamma_1* \gamma_2)* \gamma_3=\gamma_1*(\gamma_2* \gamma_3)\]
as soon as the composite exists.

\bd \label{composition_map} The map $\gamma_1*\gamma_2$ is called the {\rm composition} of $\gamma_1$ and $\gamma_2$. The composite 
\[\xymatrix@C=3em{\gamma_1 *_N \gamma_2: [0,1] \fr{N:t\mapsto 2t}& [0,2]
\fr{\gamma_1*\gamma_2}& U}\] 
is called the {\rm normalized composition}.  
\ed

\bd \cite{mdtop} A {\rm multipointed $d$-space $X$} is a triple
$(|X|,X^0,\P^{\mathcal{G}}X)$ where
\begin{itemize}
\item The pair $(|X|,X^0)$ is a multipointed space.
\item The set $\P^{\mathcal{G}}X$ is a set of continous maps from $[0,1]$ to $|X|$ called the {\rm execution paths}, satisfying the following axioms:
\begin{itemize}
\item For any execution path $\gamma$, one has $\gamma(0),\gamma(1)\in X^0$.
\item Let $\gamma$ be an execution path of $X$. Then any composite $\gamma.\phi$ with $\phi\in \mathcal{G}$ is an execution path of $X$.
\item Let $\gamma_1$ and $\gamma_2$ be two composable execution paths of $X$; then the normalized composition $\gamma_1 *_N \gamma_2$ is an execution path of $X$.
\end{itemize}
\end{itemize}
A map $f:X\to Y$ of multipointed $d$-spaces is a map of multipointed spaces from $(|X|,X^0)$ to $(|Y|,Y^0)$ such that for any execution path $\gamma$ of $X$, the map $f. \gamma$ is an execution path of $Y$.  The category of multipointed $d$-spaces is denoted by $\ptop{\mathcal{G}}$. The subset of execution paths from $\alpha$ to $\beta$ is the set of $\gamma\in\P^{\mathcal{G}} X$  such that $\gamma(0)=\alpha$ and $\gamma(1)=\beta$; it is denoted by $\P^{\mathcal{G}}_{\alpha,\beta} X$. It is equipped with the kelleyfication of the initial topology making the inclusion $\P^{\mathcal{G}}_{\alpha,\beta} X\subset \ttop([0,1],|X|)$ is continuous. \ed

\bd
Let $X$ be a multipointed $d$-space $X$. Let $\P^{\mathcal{G}} X$ be the topological space \[\P^{\mathcal{G}} X = \bigsqcup_{(\alpha,\beta)\in X^0\p X^0} \P_{\alpha,\beta}^{\mathcal{G}}X.\]
\ed

The category of multipointed $d$-spaces $\ptop{\mathcal{G}}$ is locally presentable and the forgetful functor $X\mapsto \omega(|X|)$ is topological and fibre-small by \cite[Theorem~3.5]{mdtop}. The following examples play an important role in the sequel. 
\begin{enumerate}[leftmargin=*]
\item Any set $E$ will be identified with the multipointed $d$-space $(E,E,\varnothing)$.
\item The \textit{topological globe of $Z$}, which is denoted by $\glob^{\mathcal{G}}(Z)$, is the multipointed $d$-space defined as follows
\begin{itemize}
\item the underlying topological space is the quotient space \[\frac{\{\widehat{0},\widehat{1}\}\sqcup (Z\p[0,1])}{(z,0)=(z',0)=\widehat{0},(z,1)=(z',1)=\widehat{1}}\]
\item the set of states is $\{\widehat{0},\widehat{1}\}$
\item the set of execution paths is the set of continuous maps \[\{t\mapsto (x,\phi(t))\mid t\in [0,1], \phi\in \mathcal{G},x\in  Z\}.\]
\end{itemize}
In particular, $\glob^{\mathcal{G}}(\varnothing)$ is the multipointed $d$-space $\{\widehat{0},\widehat{1}\} = (\{\widehat{0},\widehat{1}\},\{\widehat{0},\widehat{1}\},\varnothing)$. 
\item The \textit{directed segment} is the multipointed $d$-space $\vI^{\mathcal{G}}=\glob^{\mathcal{G}}(\{0\})$. 
\item The multipointed $d$-space $\vec{[\ell_1,\ell_2]}$ where $\ell_1<\ell_2$ are two real numbers has the underlying space the segment $[\ell_1,\ell_2]$, the set of states $\{\ell_1,\ell_2\}$ and the unique space of execution paths $\P_{\ell_1,\ell_2}^{\mathcal{G}}\vec{[\ell_1,\ell_2]} = \{[0,1]\iso^+ [\ell_1,\ell_2]\}$. 
\end{enumerate}

\bp \label{overmtop} The mapping $\Omega:X\mapsto (|X|,X^0)$ induces a functor from $\ptop{\mathcal{G}}$ to $\mtop$ which is topological and fibre-small. \ep

\bpf The statement is very close to the statement of \cite[Proposition~3.6]{mdtop}. The proof of the latter proposition uses the final structure. We prefer to use the $\Omega$-initial structure because it will be reused in Corollary~\ref{cstr-cocyl}. Let $(|X|,X^0)$ be a multipointed space. Consider a cone (which can be large) $(f_i:(|X|,X^0)\to  \Omega(X_i))_{i\in I}$. For all $(\alpha,\beta)\in X^0\p X^0$, consider the set of paths \[P_{\alpha,\beta}=\{\gamma\in \top([0,1],|X|)\mid \gamma(0),\gamma(1)\in X^0 \hbox{ and }\forall i,f_i.\gamma \in \P^{\mathcal{G}}_{f_i(\alpha),f_i(\beta)}X_i\}.\] Let $\gamma\in P_{\alpha,\beta}$. Let $\phi\in \mathcal{G}$. Then $\gamma(\phi(0))=\gamma(0)$, $\gamma(\phi(1))=\gamma(1)$ and $f_i.\gamma.\phi\in \P^{\mathcal{G}}_{f_i(\alpha),f_i(\beta)}X_i$ for all $i$ by definition of $P_{\alpha,\beta}$. It also means that $\gamma.\phi\in P_{\alpha,\beta}$. Let $\gamma_1\in P_{\alpha,\alpha'}$ and $\gamma_2\in P_{\alpha',\alpha''}$. Then $f_i.(\gamma_1 *_N \gamma_2)=(f_i.\gamma_1) *_N (f_i.\gamma_2)$ for all $i$ by definition of $*_N$. Therefore  $f_i.(\gamma_1 *_N \gamma_2)\in \P^{\mathcal{G}}_{f_i(\alpha),f_i(\alpha'')}X_i$ for all $i$. We deduce that $\gamma_1 *_N \gamma_2\in P_{\alpha,\alpha''}$ by definition of $P_{\alpha,\alpha''}$. We deduce that the family of $(P_{\alpha,\beta})$ yields a structure of multipointed $d$-space on $(|X|,X^0)$ and it is clearly the biggest one because all $f_i$ must be lifted to maps of multipointed $d$-spaces. It is therefore the $\Omega$-initial structure.
\epf

\begin{nota}
Let $u\in [0,1]$. Let $(Z,Z^0)\in\mtop$. Let $\ev_u:(\ttop([0,1],Z),Z^0)\to (Z,Z^0)$ be the evaluation at $u$ where $Z^0$ is identified to the corresponding set of constant maps of $\ttop([0,1],Z)$.
\end{nota}

\begin{cor} \label{cstr-cocyl} Let $X$ be a multipointed $d$-space. Let $\cocyl^{\mathcal{G}}(X)$ be the $\Omega$-initial lift of the cone $(\ev_u:(\ttop([0,1],|X|),X^0)\to (\Omega(X))_{u\in [0,1]})$ where $X^0$ is identified to the corresponding the set of constant maps of $\ttop([0,1],|X|)$. Then the space of execution paths of $\cocyl^{\mathcal{G}}(X)$ from $\alpha$ to $\beta$ is the topological space $\ttop([0,1],\P^{\mathcal{G}}_{\alpha,\beta}X)$. 
\end{cor}

\bpf
By construction of the $\Omega$-initial structure explained in the proof of Proposition~\ref{overmtop}, we have the equality of sets
\[\P^{\mathcal{G}}_{\alpha,\beta}\cocyl^{\mathcal{G}}(X)=\{\gamma\in \top([0,1],\ttop([0,1],|X|))\mid \forall u\in [0,1],\ev_u.\gamma\in \P^{\mathcal{G}}_{\alpha,\beta}X\}.\]
By endowing the two members of the equality with their topology (i.e. their initial structure in $\top$ making the inclusion into $\ttop([0,1],\ttop([0,1],|X|))$ continuous), we obtain the homeomorphism 
\[\P^{\mathcal{G}}_{\alpha,\beta}\cocyl^{\mathcal{G}}(X)\iso\{\gamma\in \ttop([0,1],\ttop([0,1],|X|))\mid \forall u\in [0,1],\ev_u.\gamma\in \P^{\mathcal{G}}_{\alpha,\beta}X\}.\]
The point is that $\top$ is cartesian closed. Therefore, we can switch the left-hand copy and the right-hand copy of the segment $[0,1]$ in $\ttop([0,1],\ttop([0,1],|X|))$. We obtain the homeomorphism  $\P^{\mathcal{G}}_{\alpha,\beta}\cocyl^{\mathcal{G}}(X) \iso \ttop([0,1],\P^{\mathcal{G}}_{\alpha,\beta}X)$ by taking $\gamma$ to the continuous mapping $u\mapsto \ev_u.\gamma$. 
\epf

\begin{nota}
The map of multipointed $d$-spaces induced by \[\ev_u:(\ttop([0,1],|X|),X^0)\to (|X|,X^0)\] is denoted by $\pi_u:\cocyl^{\mathcal{G}}(X)\to X$.
\end{nota}

\bp Let $X$ be a multipointed $d$-space. There exists a unique map $\tau : X \to \cocyl^{\mathcal{G}}(X)$ of multipointed $d$-spaces such thay the underlying map of multipointed spaces $\tau:(|X|,X^0)\to (\ttop([0,1],|X|),X^0)$ takes $x\in |X|$ to the constant path $\tau(x):t\mapsto x$ of $X$.
\ep

\bpf 
Since the functor $\Omega:X\mapsto (|X|,X^0)$ is topological, there exists at most one such a map. Let $\gamma:[0,1]\to |X|$ be an element of $\P^{\mathcal{G}}_{\alpha,\beta}X$. Then for all $t,u\in [0,1]$, one has $ev_u.\tau(\gamma(t))=\gamma(t)$. It means that $\ev_u.\tau.\gamma=\gamma \in \P^{\mathcal{G}}_{\alpha,\beta}X$. 
\epf

\begin{cor}
The mapping $X\mapsto \cocyl^{\mathcal{G}}(X)$ gives rise to a well-defined functor from $\ptop{\mathcal{G}}$ to itself. The map $\pi=(\pi_0,\pi_1)$ together with the map $\tau:X\to\cocyl^{\mathcal{G}}(X)$ above defined gives rise to a path functor, i.e. the composite $(\pi_0,\pi_1).\tau$ is the codiagonal.
\end{cor}

\begin{nota} Let $X$ be a multipointed $d$-space. Then the pair $(X^0,(\P^{\mathcal{G}}_{\alpha,\beta}X)_{\alpha,\beta})$ is a well-defined topological graph denoted by $\pf^{\mathcal{G}}(X)$. \end{nota}

\bp \label{preprediag} Let $U$ be a topological space. Let $X$ be a multipointed $d$-space. Then we have the natural bijection 
\[\ptop{\mathcal{G}}(\glob^{\mathcal{G}}(U),X) \iso \bigsqcup_{(\alpha,\beta)\in X^0\p X^0} \top(U,\P^{\mathcal{G}}_{\alpha,\beta}X).\]
\ep

\bpf A map of multipointed $d$-spaces from $\glob^{\mathcal{G}}(U)$ to $X$ is characterized by the choice of two states $\alpha$ and $\beta$ of $X$ for the image of $\widehat{0}$ and $\widehat{1}$ respectively and by a continuous map $f$ from $|\glob^{\mathcal{G}}(U)|$ to $X$ such that $f(u,-)\in \P^{\mathcal{G}}_{\alpha,\beta}X$ for all $u\in [0,1]$. In other terms, the mapping $f\mapsto (u\mapsto f(u,-))$ yields a natural set map 
\[\ptop{\mathcal{G}}(\glob^{\mathcal{G}}(U),X) \longrightarrow \bigsqcup_{(\alpha,\beta)\in X^0\p X^0} \top(U,\P^{\mathcal{G}}_{\alpha,\beta}X).\]
Conversely, consider an element $g\in \top(U,\P^{\mathcal{G}}_{\alpha,\beta}X)$ for some $(\alpha,\beta)\in X^0\p X^0$. Then the mapping $(t,u)\mapsto g(u)(t)$ induces a map of multipointed $d$-spaces from $\glob^{\mathcal{G}}(U)$ to $X$. The proof is complete because $\top$ is cartesian closed. 
\epf

\bp The mapping $X\mapsto\pf^{\mathcal{G}}(X)$ induces a well-defined functor from $\ptop{\mathcal{G}}$ to $\predtop$. It is a right adjoint. \ep

\bpf Roughly, the left adjoint is the free multipointed $d$-space generated by a topological graph. The left adjoint $\pf^{\mathcal{G}}_!:\predtop\to\ptop{\mathcal{G}}$ is constructed as follows. Let $X=(X^0,(X_{\alpha,\beta}))$ be a topological graph. We start from the set $X^0$ equipped with the discrete topology. We add a topological globe $\glob^{\mathcal{G}}(X_{\alpha,\beta})$ with $\widehat{0}$ identified with $\alpha$ and $\widehat{1}$ identified with $\beta$ for each $(\alpha,\beta)\in X^0\p X^0$. We obtain a multipointed $d$-space $\pf^{\mathcal{G}}_!(X)$. A map $f$ of multipointed $d$-spaces from $\pf^{\mathcal{G}}_!(X)$ to $Y$ is equivalent to choosing a set map from $\pf^{\mathcal{G}}_!(X)^0=X^0$ to $Y^0$ and for each $(\alpha,\beta)\in X^0\p X^0$ a map of multipointed $d$-spaces from $\glob^{\mathcal{G}}(X_{\alpha,\beta})$ to $Y$, which is equivalent by Proposition~\ref{preprediag} to choosing a map from $X_{\alpha,\beta}$ to $\P_{f(\alpha),f(\beta)}^{\mathcal{G}}Y$. 
\epf

\bth \label{three} Let $(\C,\W,\F)$ be one of the three model structures \[(\C_q,\W_q,\F_q), (\C_{\overline{h}},\W_h,\F_h),(\C_m,\W_m,\F_m)\] of $\top$. Then there exists a unique model structure on $\ptop{\mathcal{G}}$ such that: 
\begin{itemize}
\item A map of multipointed $d$-spaces $f:X\to Y$ is a weak equivalence if and only if $f^0:X^0\to Y^0$ is a bijection and for all $(\alpha,\beta)\in X^0\p X^0$, the continuous map $\P_{\alpha,\beta}^{\mathcal{G}}X\to \P_{f(\alpha),f(\beta)}^{\mathcal{G}}X$ belongs to $\W$.
\item A map of multipointed $d$-spaces $f:X\to Y$ is a fibration if and only if for all $(\alpha,\beta)\in X^0\p X^0$, the continuous map $\P_{\alpha,\beta}^{\mathcal{G}}X\to \P_{f(\alpha),f(\beta)}^{\mathcal{G}}X$ belongs to $\F$.
\end{itemize}
Moreover, this model structure is accessible and all objects are fibrant.
\eth

\bpf 
A model structure is characterized by its fibrations and its weak equivalences. For all topological spaces $U$, the constant path map $\tau:U\to \ttop([0,1],U)$ is equal to the composite $U\iso \ttop(\{0\},U)\to \ttop([0,1],U)$, and is therefore a homotopy equivalence. By Proposition~\ref{preparation-util}, the map $\pi=(\pi_0,\pi_1):\ttop([0,1],U)\to \ttop(\{0,1\},U)\iso U\p U$ (the latter homeomorphism coming from the fact that $\top$ is cartesian closed) is a fibration in the three model structures. We deduce that for all multipointed $d$-spaces $X$ and all $(\alpha,\beta)\in X^0\p X^0$, the continuous map $\tau:\P^{\mathcal{G}}_{\alpha,\beta}X \to \ttop([0,1],\P^{\mathcal{G}}_{\alpha,\beta}X)$ belongs to $\W$ and the continuous map $\pi:\ttop([0,1],\P^{\mathcal{G}}_{\alpha,\beta}X)\to \P^{\mathcal{G}}_{\alpha,\beta}X\p \P^{\mathcal{G}}_{\alpha,\beta}X$ belongs to $\F$. By Corollary~\ref{cstr-cocyl}, we deduce that the factorization of the diagonal \[\xymatrix@C=3em{X \fr{\tau}& \cocyl^{\mathcal{G}}(X) \fr{\pi}& X\p X}\] satisfies the hypotheses of Theorem~\ref{GHKRS} applied to the right adjoint $\pf^{\mathcal{G}}:\ptop{\mathcal{G}}\to \predtop$. The proof is complete thanks to Corollary~\ref{qhm-graph}.
\epf

\bd The three model structures on $\ptop{\mathcal{G}}$ are called the {\rm q-model structure}, the {\rm h-model structure} and the {\rm m-model structure} respectively. They are denoted by $(\ptop{\mathcal{G}})_q$, $(\ptop{\mathcal{G}})_h$ and $(\ptop{\mathcal{G}})_m$ respectively. 
\ed

\bth \label{leftdetspace} The q-model structure of $\ptop{\mathcal{G}}$ is combinatorial and left determined. It coincides with the combinatorial model structure of \cite{mdtop}. A set of generating cofibrations is $\{\glob^{\mathcal{G}}(\mathbf{S}^{n-1}) \subset \glob^{\mathcal{G}}(\mathbf{D}^{n})\mid n\geq 0\} \cup \{C:\varnothing \to \{0\},R:\{0,1\} \to \{0\}\}$. 
\eth

\bpf The q-model structure of $\ptop{\mathcal{G}}$ coincides with the model structure of \cite{mdtop} since fibrations and weak equivalences determine a model structure. Therefore it is combinatorial. It is left determined by a proof similar to the one of \cite[Theorem~4.3]{leftdetflow} for the category of flows: it suffices to replace $\{[0,1],Y\}_S$ (which is denoted by $\cocyl(Y)$ in the proof of Theorem~\ref{three2}) by $\cocyl^{\mathcal{G}}(Y)$, $\glob$ by $\glob^{\mathcal{G}}$ and $\P Y$ by $\P^{\mathcal{G}}Y$.
\epf

\bth \label{mdtop-mixed} The m-model structure of $\ptop{\mathcal{G}}$ is the mixed model structure of the q-model structure and the h-model structure in the sense of \cite[Theorem~2.1]{mixed-cole}. \eth

\bpf A model structure is characterized by its class of weak equivalences and by its class of fibrations. The m-model 
structure of $\ptop{\mathcal{G}}$ is therefore the unique model structure such that a map of multipointed $d$-spaces $f:X\to Y$ is
\begin{itemize}
\item a weak equivalence if and only if it is a weak equivalence of the q-model structure of $\ptop{\mathcal{G}}$.
\item a fibration if and only if it is a fibration of the h-model structure of $\ptop{\mathcal{G}}$.
\end{itemize}
Hence the proof is complete. 
\epf

\bp \label{implication-mdtop} There are the implications $\hbox{q-cofibrant} \Rightarrow \hbox{m-cofibrant} \Rightarrow \hbox{h-cofibrant}$ for $\ptop{\mathcal{G}}$. The identity functor yields a Quillen equivalence \[\id:(\ptop{\mathcal{G}})_{q} \dashv (\ptop{\mathcal{G}})_{m}:\id.\]
\ep

\bpf
The first assertion is a consequence of \cite[Corollary~3.7]{mixed-cole}. The second assertion is obvious. 
\epf

\bp \label{mdtop-notcofibrant} There exists a multipointed $d$-space which is not cofibrant in any of the three model structures of Theorem~\ref{three}.
\ep

\begin{figure}
\[
\xymatrix@R=1em@C=3em
{
&& \stackrel{(0,a)}{\bullet} \ar@{--}[dddd] \ar@/^10pt/@{->}[rrd]^-{} & \\
\overline{X} \ar@{->}[ddd]_-{p} &\stackrel{\widehat{0}}{\bullet} \ar@{->}[rrr]\ar@/^5pt/@{->}[rrr]\ar@/_5pt/@{->}[rrr]
\ar@/^10pt/@{->}[rrr]\ar@/_10pt/@{->}[rrr]
\ar@/^15pt/@{->}[rrr]\ar@/_15pt/@{->}[rrr]
\ar@/^20pt/@{->}[rrr]\ar@/_20pt/@{->}[rrr]
\ar@/^25pt/@{->}[rrr]\ar@/_25pt/@{->}[rrr]
\ar@{--}[ddd]\ar@/^5pt/@{->}[ru]^-{} \ar@/_10pt/@{->}[rrd]_-{} &&& \stackrel{\widehat{1}}{\bullet} \ar@{--}[ddd]\\
&&& \stackrel{(1,b)}{\bullet} \ar@{--}[dd]\ar@/_5pt/@{->}[ru]_-{}& \\
\\
X &\stackrel{0}{\bullet} \ar@{-}[r] & \stackrel{a}{\bullet} \ar@{-}[r] & \stackrel{b}{\bullet} \ar@{->}[r] & \stackrel{1}{\bullet}
}
\]
\caption{Symbolic representation of $p:\overline{X}\to X$}
\label{Z}
\end{figure}
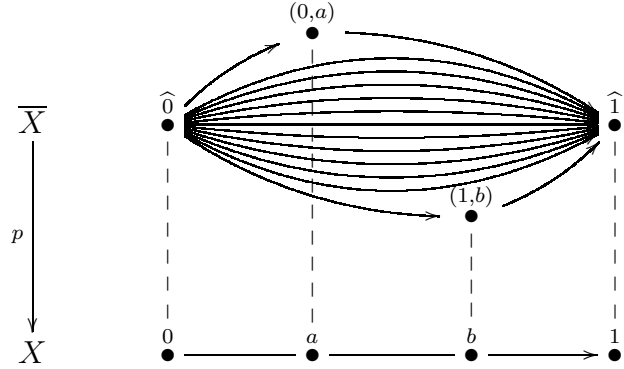

\bpf
Consider the poset $\widehat{P}$ consisting of the set $\{0,a,b,1\}$ equipped with the partial order $0<a<1$ and $0<b<1$: $a$ and $b$ are not comparable. Let $X$ be the multipointed $d$-space defined as follows. Let $|X|=[0,1]$. Let $X^0=\{0,a,b,1\}$ with $a=1/3$ and $b=2/3$. Let $\P_{\alpha,\beta}^{\mathcal{G}}X = \{[0,1] \iso^+ [\alpha,\beta]\}$ if and only if $\alpha<\beta$ in $\widehat{P}$ and  $\P_{\alpha,\beta}^{\mathcal{G}}X =\varnothing$ otherwise. These data clearly satisfy the axioms of multipointed $d$-space. Consider the multipointed $d$-space $\overline{X}$ defined as follow: 
\begin{enumerate}[leftmargin=*]
\item We start from the multipointed $d$-space \[(\vec{[0,a]}*\vec{[a,1]})\sqcup (\vec{[0,b]}*\vec{[b,1]})\] where the symbol $*$ means that the two copies of $a$ (of $b$ resp.) are identified.
\item We make the identifications $0=0$ and $1=1$, we obtain a multipointed $d$-space $Z$ whose underlying space is homeomorphic to $\mathbf{S}^1$. We consider the pushout diagram of multipointed $d$-spaces
\[
\xymatrix
{
\glob^{\mathcal{G}}(\{0,1\})\fd{} \fr{\phi} & Z \fd{} \\
\glob^{\mathcal{G}}([0,1]) \fr{} & \cocartesien \overline{X}
}
\]
such that $\phi(\widehat{0})=0$, $\phi(\widehat{1})=1$, $\phi$ maps $(0,t)$ to $t\in|\vec{[0,a]}*\vec{[a,1]}|$ and maps $(1,t)$ to $t\in|\vec{[0,b]}*\vec{[b,1]}|$. Intuitively, the hole in the middle of $Z$ is filled by a homotopy between the execution path $(0,1)$ of $\vec{[0,a]}*\vec{[a,1]}$ and the execution path $(0,1)$ of $\vec{[0,b]}*\vec{[b,1]}$. It is depicted in Figure~\ref{Z}. 
\end{enumerate}
The projection map $(z,t)\to t$ from $[0,1]\p [0,1]$ to $[0,1]$ induces a map of multipointed $d$-spaces $p:\overline{X}\to X$ preserving the set of states. It is depicted in Figure~\ref{Z} as well. The maps $\P_{\alpha,\beta}^{\mathcal{G}}\overline{X} \to \P_{p(\alpha),p(\beta)}^{\mathcal{G}}X$ are either $\id_\varnothing$, $\id_{\mathcal{G}}$, and for $(\alpha,\beta)=(\widehat{0},\widehat{1})$, it is the projection map $[0,1]\p \mathcal{G} \to \mathcal{G}$ which is homotopy equivalence and a h-fibration of $\top$. Therefore the map $p:\overline{X}\to X$ is a trivial fibration of the h-model structure of $\ptop{\mathcal{G}}$. 

If $X$ was h-cofibrant, then there would exist a section $s:X\to \overline{X}$ of $p:\overline{X}\to X$. Since $p$ induces a bijection from $\overline{X}^0$ to $X^0$, the map of multipointed $d$-spaces $s:X\to \overline{X}$ must induce a bijection from $X^0$ to $\overline{X}^0$. It means that $s(0)=\widehat{0}$, $s(a)=(0,a)$, $s(b)=(1,b)$ and $s(1)=\widehat{1}$. Let $\gamma \in \P_{a,1}^{\mathcal{G}}X$. Then the continuous map $\gamma:[0,1] \to [a,1]$ is a nondecreasing homeomorphism. Since $s:X\to \overline{X}$ is a map of multipointed $d$-spaces, the composite $s.\gamma$ belongs to $\P_{(0,a),\widehat{1}}^{\mathcal{G}}\overline{X}$. The point is that $s(\gamma(\gamma^{-1}(b))) = s(b) = (1,b)$. It means that $s.\gamma$ is an execution path of $\overline{X}$ from $(0,a)$ to $\widehat{1}$ passing by $(1,b)$. Such an execution path does not exist in $\overline{X}$ by construction: all execution paths inside the globe are parallel to the boundary indeed. Contradiction. We deduce that $X$ is not h-cofibrant. The proof is complete thanks to Proposition~\ref{implication-mdtop}. 
\epf

\section{Flow}
\label{homotopy-flow}

\bd \cite{model3} A {\rm flow} $X$ consists of a topological space $\P X$ of execution paths, a discrete space $X^0$ of states, two continuous maps $s$ and $t$ from $\P X$ to $X^0$ called the source and target map respectively, and a continuous and associative map \[*:\{(x,y)\in \P X\p \P X; t(x)=s(y)\}\longrightarrow \P X\] such that $s(x*y)=s(x)$ and $t(x*y)=t(y)$.  A morphism of flows $f:X\longrightarrow Y$ consists of a set map $f^0:X^0\longrightarrow Y^0$ together with a continuous map $\P f:\P X\longrightarrow \P Y$ such that $f(s(x))=s(f(x))$, $f(t(x))=t(f(x))$ and $f(x*y)=f(x)*f(y)$. The corresponding category is denoted by $\dtop$. Let $\P_{\alpha,\beta}X = \{x\in \P X\mid s(x)=\alpha \hbox{ and } t(x)=\beta\}$.
\ed

The category $\dtop$ is locally presentable by \cite[Theorem~7.7]{mdtop}. Three examples of flows are important for the sequel: 
\begin{enumerate}[leftmargin=*]
\item For a topological space $X$, let $\glob(X)$ be the flow defined by $\glob(X)^0=\{0,1\}$ and $\P \glob(X)=X$ with $s=0$ and $t=1$. This flow has no composition law.
\item The flow $\vI$ is by definition $\glob(\{0\})$.
\item Let $(P,\leq)$ be a poset. Then it can be viewed as a flow with $\P_{\alpha,\beta}P$ equal to the singleton $\{(\alpha,\beta)\}$ if and only if $\alpha<\beta$ and empty otherwise. In particular, a set can be viewed as a flow without execution paths.
\end{enumerate}

\begin{nota} Let $X$ be a flow. Then the pair $(X^0,(\P_{\alpha,\beta}X)_{\alpha,\beta})$ is a well-defined topological graph denoted by $\pf(X)$. \end{nota}

\bp The mapping $X\mapsto\pf(X)$ induces a well-defined functor from $\dtop$ to $\predtop$. It is a right adjoint. \ep

\bpf Roughly, the left adjoint is the free flow generated by a topological graph. The left adjoint $\pf_!:\predtop\to\dtop$ is constructed as follows. Let $X=(X^0,(X_{\alpha,\beta}))$ be a topological graph. The set of states of $\pf_!(X)$ is $X^0$. For $\alpha,\beta\in X^0$, let 
\[\P_{\alpha,\beta}X = \bigsqcup_{\tiny\begin{array}{c}
(\alpha_1,\dots,\alpha_n)\in (X^0)^{n}\\
n\geq 2,\alpha_1=\alpha,\alpha_n=\beta
\end{array}} X_{\alpha_1,\alpha_2} \p \dots \p X_{\alpha_{n-1},\alpha_n}.\]
The composition law is defined by concatening tuples: \[(x_1,\dots,x_m)*(y_1,\dots,y_n)=(x_1,\dots,x_m,y_1,\dots,y_n)\]
We obtain a flow $\pf_!(X)$. A map $f$ of flows from $\pf_!(X)$ to $Y$ is equivalent to choosing a set map from $\pf_!(X)^0=X^0$ to $Y^0$ and for each $(\alpha,\beta)\in X^0\p X^0$ a continous map from $X_{\alpha,\beta}$ to $Y_{f(\alpha),f(\beta)}$. 
\epf

\bth \label{three2} Let $(\C,\W,\F)$ be one of the three model structures \[(\C_q,\W_q,\F_q), (\C_{\overline{h}},\W_h,\F_h),(\C_m,\W_m,\F_m)\] of $\top$. Then there exists a unique model structure on $\dtop$ such that: 
\begin{itemize}
\item A map of flows $f:X\to Y$ is a weak equivalence if and only if $f^0:X^0\to Y^0$ is a bijection and for all $(\alpha,\beta)\in X^0\p X^0$, the continuous map $\P_{\alpha,\beta}X\to \P_{f(\alpha),f(\beta)}X$ belongs to $\W$.
\item A map of multipointed $d$-spaces $f:X\to Y$ is a fibration if and only if for all $(\alpha,\beta)\in X^0\p X^0$, the continuous map $\P_{\alpha,\beta}X\to \P_{f(\alpha),f(\beta)}X$ belongs to $\F$.
\end{itemize}
Moreover, this model structure is accessible and all objects are fibrant.
\eth

\bpf[Sketch of proof]
The proof is similar to the proof of Theorem~\ref{three}. Roughly speaking, it suffices to replace everywhere $\P_{\alpha,\beta}^{\mathcal{G}}X$ by $\P_{\alpha,\beta}X$ and to use the right adjoint $\pf:\dtop\to\predtop$. We also have to use the path functor $\cocyl:\dtop\to\dtop$ defined on objects by $\cocyl(X)^0:=X^0$, for all $(\alpha,\beta)\in X^0\p X^0$, $\P_{\alpha,\beta}\cocyl(X) := \ttop([0,1],\P_{\alpha,\beta}X)$ with an obvious definition of the composition law. It is the flow denoted by $\{[0,1],X\}_S$ in \cite[Notation~7.6]{model3} and in \cite[Notation~3.8]{leftdetflow}. 
\epf

\bd The three model structures on $\dtop$ are called the {\rm q-model structure}, the {\rm h-model structure} and the {\rm m-model structure} respectively. They are denoted by $(\dtop)_q$, $(\dtop)_h$ and $(\dtop)_m$ respectively.
\ed

\bth \label{leftdetflow} The q-model structure of $\dtop$ is combinatorial and left determined. It coincides with the combinatorial model structure of \cite{model3}. A set of generating cofibrations is $\{\glob(\mathbf{S}^{n-1}) \subset \glob(\mathbf{D}^{n})\mid n\geq 0\} \cup \{C:\varnothing \to \{0\},R:\{0,1\} \to \{0\}\}$
\eth

\bpf The q-model structure of $\dtop$ coincides with the model structure of \cite{model3} since fibrations and weak equivalences determine a model structure. Therefore it is combinatorial. It is left determined by \cite[Theorem~4.3]{leftdetflow}.
\epf

\bth \label{mflow-mixed} The m-model structure of $\dtop$ is the mixed model structure of the q-model structure and the h-model structure in the sense of \cite[Theorem~2.1]{mixed-cole}. \eth

\bpf A model structure is characterized by its class of weak equivalences and by its class of fibrations. The m-model 
structure of $\dtop$ is therefore the unique model structure such that a map of flows $f:X\to Y$ is
\begin{itemize}
\item a weak equivalence if and only if it is a weak equivalence of the q-model structure of $\dtop$.
\item a fibration if and only if it is a fibration of the h-model structure of $\dtop$.
\end{itemize}
Hence the proof is complete. 
\epf

\bp \label{impl-dtop} There are the implications $\hbox{q-cofibrant} \Rightarrow \hbox{m-cofibrant} \Rightarrow \hbox{h-cofibrant}$ for $\dtop$. The identity functor yields a Quillen equivalence \[\id:(\dtop)_{q} \dashv (\dtop)_{m}:\id.\] 
\ep

\bpf
The first assertion is a consequence of \cite[Corollary~3.7]{mixed-cole}. The second assertion is obvious. 
\epf  

\bp \label{flow-notcofibrant} There exists a flow which is not cofibrant in any of the three model structures of Theorem~\ref{three2}.
\ep

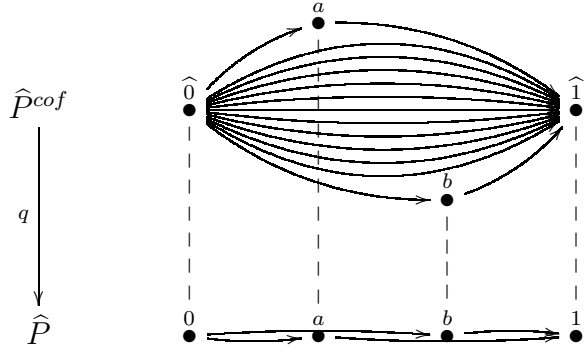
\begin{figure}
\[
\xymatrix@R=1em@C=3em
{
&& \stackrel{a}{\bullet} \ar@{--}[dddd] \ar@/^10pt/@{->}[rrd]^-{} & \\
\widehat{P}^{cof} \ar@{->}[ddd]_-{q} &\stackrel{\widehat{0}}{\bullet} \ar@{->}[rrr]\ar@/^5pt/@{->}[rrr]\ar@/_5pt/@{->}[rrr]
\ar@/^10pt/@{->}[rrr]\ar@/_10pt/@{->}[rrr]
\ar@/^15pt/@{->}[rrr]\ar@/_15pt/@{->}[rrr]
\ar@/^20pt/@{->}[rrr]\ar@/_20pt/@{->}[rrr]
\ar@/^25pt/@{->}[rrr]\ar@/_25pt/@{->}[rrr]
\ar@{--}[ddd]\ar@/^5pt/@{->}[ru]^-{} \ar@/_10pt/@{->}[rrd]_-{} &&& \stackrel{\widehat{1}}{\bullet} \ar@{--}[ddd]\\
&&& \stackrel{b}{\bullet} \ar@{--}[dd]\ar@/_5pt/@{->}[ru]_-{}& \\
\\
\widehat{P} &\stackrel{0}{\bullet} \ar@/^2pt/@{->}[rr] \ar@/_2pt/@{->}[r]& \stackrel{a}{\bullet}\ar@/_2pt/@{->}[rr]  & \stackrel{b}{\bullet} \ar@/^2pt/@{->}[r] & \stackrel{1}{\bullet}
}
\]
\caption{Symbolic representation of $q:\widehat{P}^{cof}\to \widehat{P}$}
\label{ZZ}
\end{figure}

\bpf
As in the proof of Proposition~\ref{mdtop-notcofibrant}, consider the poset $\widehat{P}$ consisting of the set $\{0,a,b,1\}$ equipped with the partial order $0<a<1$ and $0<b<1$: $a$ and $b$ are not comparable. We denoted in the same way the flow associated with the poset $\widehat{P}$. We consider a q-cofibrant replacement $\widehat{P}^{cof}$ of $\widehat{P}$ constructed as follows: 
\begin{enumerate}[leftmargin=*]
\item We start from the flow \[(\vI*\vI)\sqcup (\vI*\vI)\] where the symbol $*$ means that the final state of the left-hand flow is identified with the initial state of the right-hand flow. The middle state of the left-hand copy of $\vI*\vI$ is denoted by $a$ and the middle state of the right-hand copy of $\vI*\vI$ is denoted by $b$.
\item We make the identifications $\widehat{0}=\widehat{0}$ and $\widehat{1}=\widehat{1}$, we obtain a flow $T$. We consider the pushout diagram of multipointed $d$-spaces
\[
\xymatrix
{
\glob(\{0,1\})\fd{} \fr{\phi} & T \fd{} \\
\glob([0,1]) \fr{} & \cocartesien \widehat{P}^{cof}
}
\]
such that $\phi(\widehat{0})=\widehat{0}$, $\phi(\widehat{1})=\widehat{1}$, $\phi$ maps the path $0$ to the unique execution path of the left-hand copy of $\vI*\vI$ from $\widehat{0}$ to $\widehat{1}$ and the path $1$ to the unique execution path of the right-hand copy of $\vI*\vI$ from $\widehat{0}$ to $\widehat{1}$. It is depicted in Figure~\ref{ZZ}. 
\end{enumerate}
There is a unique map of flows $q:\widehat{P}^{cof}\to \widehat{P}$. It preserves the set of states. It is a h-fibration of $\dtop$ since all topological spaces $\P_{\alpha,\beta}\widehat{P}$ are either singleton, or empty and in this case $\P_{\alpha,\beta}\widehat{P}^{cof}$ is (and must be) empty as well. It is a trivial fibration of $(\dtop)_h$ because all nonempty path spaces of $\widehat{P}^{cof}$ are contractible. 

If $\widehat{P}$ was h-cofibrant, then there would exist a section $s:\widehat{P}\to \widehat{P}^{cof}$ of $q$. Since $q:\widehat{P}^{cof}\to \widehat{P}$ induces a bijection between the set of states,  we would have $s(0)=\widehat{0}$, $s(a)=a$, $s(b)=b$ and $s(1)=\widehat{1}$. The only execution path $(0,a)$ of $P$ from $0$ to $a$ is mapped to the only execution path $s(0,a)$ of $P^{cof}$ from $\widehat{0}$ to $a$. In the same way, $s(a,1)$ is the only execution path of $\widehat{P}^{cof}$ from $a$ to $\widehat{1}$, $s(0,b)$ is the only execution path of $P^{cof}$ from $\widehat{0}$ to $b$ and finally $s(b,1)$ is the only execution path of $\widehat{P}^{cof}$ from $b$ to $\widehat{1}$. We obtain $s(0,a)*s(a,1) = s(0,1) = s(0,b)*s(b,1)$. Contradiction because there is only a homotopy in $\widehat{P}^{cof}$ between $s(0,a)*s(a,1)=\phi(0)$ and $s(0,b)*s(b,1)=\phi(1)$. The proof is complete thanks to Proposition~\ref{impl-dtop}. 
\epf 

\section{Path space functor and m-cofibrancy}
\label{pathspace}

Let us mention the erratum in \cite{leftproperflow} correcting some proofs of \cite{model2}. We conclude this paper by explaining why the m-model structures of multipointed $d$-spaces and of flows are better behaved than their q-model structures. Let us start with the following observation: 

\bth \label{flow-cof}  Let $X$ be a q-cofibrant flow. Then the space of execution paths $\P X$ is q-cofibrant. 
\eth

\bpf This fact, stated in various papers before this one, has a correct proof in \cite{leftproperflow}.
\epf

The analogue fact for multipointed $d$-spaces is wrong. Indeed, the multipointed $d$-space $\glob^{\mathcal{G}}(\mathbf{D}^{1})$ is q-cofibrant. Its space of paths is equal to $\mathbf{D}^1 \p \mathcal{G}$ which is far from being q-cofibrant in $\top$. However, it is a m-cofibrant space by \cite[Corollary~3.7]{mixed-cole} because the topological group $\mathcal{G}$ is contractible. It turns out that this phenomenon is general. We need first to recall some results of \cite{mdtop} and \cite{model2} to facilitate the reading of the proof for a reader who would not be familiar with our work.

\begin{nota} 
Let $X$ be a multipointed $d$-space. For every $(\alpha,\beta)\in X^0\p X^0$,  let $\P_{\alpha,\beta}X := \P^{\mathcal{G}}_{\alpha,\beta}X/\mathcal{G}$ be the quotient of the space $\P^{\mathcal{G}}_{\alpha,\beta}X$ by the actions of $\mathcal{G}$ equipped with the final structure, i.e. the final topology.
\end{nota}

Let $X$ be a multipointed $d$-space. Then there exists a unique flow $cat(X)$ with $cat(X)^0=X^0$, $\P_{\alpha,\beta}cat(X)= \P_{\alpha,\beta}X$ for every $(\alpha,\beta)\in X^0\p X^0$ and the composition law $*:\P_{\alpha,\beta}X \p \P_{\beta,\gamma}X \rightarrow \P_{\alpha,\gamma}X$ is for every triple $(\alpha,\beta,\gamma)\in X^0\p X^0\p X^0$ the unique map making the following diagram commutative:
\[
\xymatrix{
\P_{\alpha,\beta}^{\mathcal{G}}X \p \P_{\beta,\gamma}^{\mathcal{G}}X
\fr{*_N}\fd{} &
\P_{\alpha,\gamma}^{\mathcal{G}}X \fd{} \\
\P_{\alpha,\beta}X \p \P_{\beta,\gamma}X \fr{\exists !} &
\P_{\alpha,\gamma}X.}
\] 
The mapping $X \mapsto cat(X)$ induces a functor from $\ptop{\mathcal{G}}$ to $\dtop$ (see \cite[Section~7]{mdtop} for a complete exposition). In particular, for all topological $Z$, we have \[cat(\glob^{\mathcal{G}}(Z))=\glob(Z).\]

\begin{nota} 
Let $X,Y\in\ptop{\mathcal{G}}$. Let $\tptop{\mathcal{G}}(X,Y)$ be the set $\ptop{\mathcal{G}}(X,Y)$ equipped with the $\omega$-initial structure coming from the inclusion of sets \[\ptop{\mathcal{G}}(X,Y)\subset \mtop((|X|,X^0),(|Y|,Y^0)).\]
\end{nota}

\begin{nota} 
Let $X,Y\in\dtop$. Let $\tdtop(X,Y)$ be the set $\dtop(X,Y)$ equipped with the $\omega$-initial structure coming from the inclusion of sets \[\dtop(X,Y)\subset \set(X^0,Y^0)\p\top(\P X,\P Y),\]
with $\set(X^0,Y^0)$ equipped with the discrete topology.
\end{nota}

\bp \cite[Proposition~IV.3.1]{model2} \label{homocat}
Let $Z$ be a compact topological space. Let $U$ be a cellular object of the q-model structure of $\ptop{\mathcal{G}}$ (in \cite{model2}, such an object is called a globular complexe). Then the continuous map induced by the functor $cat:\ptop{\mathcal{G}}\to\dtop$ 
\[
cat:\tptop{\mathcal{G}}(\glob^{\mathcal{G}}(Z),U) \longrightarrow \tdtop(\glob(Z),cat(U))
\]
is a homotopy equivalence.
\ep

In fact, this proposition is a particular case of a more general theorem. In \cite[Theorem IV.3.10]{model2}, it is proved that $\glob^{\mathcal{G}}(Z)$ can be actually replaced by any cellular object $X$ of the q-model structure of $\ptop{\mathcal{G}}$, and $\glob(Z)$ must then be replaced by $cat(X)$. It is even proved in \cite[Theorem IV.3.14]{model2} that this map is a h-fibration of $\top$. The proofs of these theorems, written down within the category of weakly Hausdorff $k$-spaces, are still valid in our framework since they lie on three facts: 
\begin{enumerate}[leftmargin=*]
\item All maps of $\mathcal{G}$ are invertible: see the introduction for a short discussion about this hypothesis.
\item The underlying category of topological spaces must be bicomplete, cartesian closed and must contain all CW-complexes. 
\item The underlying category of topological spaces must be endowed with a h-model structure which is required for the homotopical part of the proofs which uses model category techniques.
\end{enumerate}

We are now able to generalize the observation above: 

\bth \label{goodbehave} Let $U$ be a m-cofibrant multipointed $d$-space. Then the space of paths $\P^{\mathcal{G}} U$ is m-cofibrant. 
\eth

\bpf
By Theorem~\ref{mdtop-mixed} and \cite[Corollary~3.7]{mixed-cole}, there exists a q-cofibrant multipointed $d$-space $V$ and a map $f:U\to V$ which is a weak equivalence of the h-model structure of $\ptop{\mathcal{G}}$. It means that $f$ induces a bijection from $U^0$ to $V^0$ and that for each $(\alpha,\beta)\in U^0 \p U^0$, the map $f:\P_{\alpha,\beta}^{\mathcal{G}}U \to \P_{f(\alpha),f(\beta)}^{\mathcal{G}}V$ is a homotopy equivalence. Therefore we can suppose without loss of generality that $U$ is q-cofibrant. Since any q-cofibrant object is a retract of a cellular one, we can suppose that $U$ is a cellular object of the q-model structure of $\ptop{\mathcal{G}}$. From a pushout diagram of multipointed $d$-spaces with $U_1$ (and therefore $U_2$) cellular
\[
\xymatrix{
\glob^{\mathcal{G}}(\mathbf{S}^{n-1}) \fd{}\fr{} & U_1 \fd{} \\
\glob^{\mathcal{G}}(\mathbf{D}^{n}) \fr{} & \cocartesien U_2,
}
\]
one obtains a pushout diagram of cellular flows 
\[
\xymatrix{
\glob(\mathbf{S}^{n-1}) \fd{}\fr{} & cat(U_1) \fd{} \\
\glob(\mathbf{D}^{n}) \fr{} & \cocartesien cat(U_2).
}
\]
This point is explained in the body of the proof of \cite[Theorem IV.3.10]{model2}. It is also easily seen that the functor $cat:\ptop{\mathcal{G}} \to \dtop$ preserves transfinite colimits of q-cofibrations between cellular objects. It is even the method used in \cite{model2} to construct the mapping $cat$. Note that the functor $cat:\ptop{\mathcal{G}} \to \dtop$ does not preserve colimits in general. Indeed, it does not have any right adjoint by \cite[Proposition~7.3]{mdtop} and being colimit-preserving and being a left adjoint are equivalent where the source and the target categories of a functor are locally presentable. 

These facts are sufficient to conclude the proof. The flow $cat(U)$ is cellular, and therefore q-cofibrant. By Theorem~\ref{flow-cof}, we deduce that the space $\P cat(U)$ is q-cofibrant. By Proposition~\ref{homocat} applied with $Z$ a singleton, the quotient map $\P^{\mathcal{G}}U \to \P cat(U)$ is a homotopy equivalence. By \cite[Corollary~3.7]{mixed-cole}, we obtain that $\P^{\mathcal{G}}U$ is a m-cofibrant space and the proof is complete. 
\epf

The same phenomenon holds for the category of flows: 

\bth \label{betterbehave} Let $U$ be a m-cofibrant flow. Then the space of paths $\P U$ is m-cofibrant. 
\eth

\bpf[Sketch of proof] There exists a map $f:U\to V$ which a weak equivalence of the h-model structure of $\dtop$ towards a q-cofibrant flow $V$. Thus $\P U$ and $\P V$ are homotopy equivalent. By Theorem~\ref{flow-cof}, the space $\P V$ is q-cofibrant.  By \cite[Corollary~3.7]{mixed-cole}, the space $\P U$ is therefore m-cofibrant. 
\epf

%\bibliographystyle{alpha} 
%\bibliography{EnrichedGph}

\end{document}